\documentclass[10pt,twoside]{amsart}
\usepackage{geometry}
\usepackage[english]{babel}
\usepackage{graphicx}
\usepackage{amsmath}
\usepackage{amsfonts}
\usepackage{amssymb}

\begin{document}
\bibliographystyle{elsarticle-num}
\numberwithin{equation}{section}
\newtheorem{teo}{Theorem}
\newtheorem{lemma}{Lemma}
\newtheorem{defi}{Definition}
\newtheorem{coro}{Corollary}
\newtheorem{prop}{Proposition}
\newtheorem{remark}{Remark}
\newtheorem{scho}{Scholium}
\newtheorem{example}{Example}
\newtheorem{open}{Open Problem}
\numberwithin{lemma}{section}
\numberwithin{prop}{section}
\numberwithin{teo}{section}
\numberwithin{defi}{section}
\numberwithin{coro}{section}
\numberwithin{figure}{section}
\numberwithin{remark}{section}
\numberwithin{scho}{section}
\numberwithin{example}{section}

\title{Geometric constants for quantifying the difference between orthogonality types}

\author{Vitor Balestro, Horst Martini, and Ralph Teixeira}
\address [V. Balestro] {CEFET/RJ Campus Nova Friburgo - Nova Friburgo - Brazil
\& Instituto de Matem\'{a}tica e Estat\'{i}stica - UFF - Niter\'{o}i - Brazil}
\email{vitorbalestro@mat.uff.br}
\address [H. Martini] {Fakult\"at f\"ur Mathematik - Technische Universit\"at Chemnitz - 09107 Chemnitz - Germany;
Dept. of Applied Mathematics - Harbin University of Science and Technology, 150080 Harbin - China}
\email{martini@mathematik.tu-chemnitz.de}
\address [R. Teixeira] {Instituto de Matem\'{a}tica e Estat\'{i}stica - UFF - Niter\'{o}i - Brazil}
\email{ralph@mat.uff.br}

\begin{abstract} This paper is devoted to introduce new geometric constants that quantify the difference between Roberts orthogonality and Birkhoff orthogonality in normed planes. We start by characterizing Roberts orthogonality in two different ways: via bisectors of two points and using certain linear transformations. Each of these characterizations yields one of those geometric constants that we will study.
\end{abstract}

\maketitle

\section{Introduction}\label{intro}

In \cite{wu} the authors introduced a geometric constant to measure the difference between Birkhoff orthogonality and isosceles orthogonality, and in \cite{Pa-Wu} analogous results for Birkhoff orthogonality and Roberts orthogonality are obtained. The main objective of this paper is to introduce two new geometric constants for
 quantifying the difference between Birkhoff orthogonality and Roberts orthogonality and, thus, continuing the investigations from \cite{Pa-Wu}. For this purpose, we present two new characterizations of Roberts orthogonality. One of them is related to segments whose bisectors contain lines, and the other one associates this type of orthogonality to certain symmetries of the unit circle.

In order to prepare these characterizations, which are given in Section \ref{roberts}, we devote Section \ref{bisectors} to the study of the geometric structure of bisectors in normed planes. The results presented in this section are not new and can be found in \cite{jahn}. We present new proofs (which are slightly more geometric in nature) that will be useful for our aim,
and we will also refer to more references dealing with geometric properties of bisectors in normed planes.

In Section \ref{cb} we introduce the constant $c_B$ using the generalized sine function defined in \cite{szostok} and studied in \cite{bmt}. In some sense this constant estimates how far the bisector of a segment is from being (or containing) a line. In Section \ref{cs} we define the constant $c_S$, which quantifies the maximum asymmetry of the unit circle regarding directions which are Birkhoff orthogonal. The reason why both these constants can be used for estimating the difference between Roberts orthogonality and Birkhoff orthogonality becomes clear already in Section \ref{roberts}.

Let us introduce some notation. Throughout the text, $(V,||\cdot||)$ will always denote a real (\emph{normed} or) \emph{Minkowski plane}, i.e., a two-dimensional vector space over $\mathbb{R}$ endowed with a norm. Its origin will be denoted by $o$, and the letters $B$ and $S$ stand, respectively, for the \textit{unit ball} $B:=\{x \in V : ||x|| \leq 1\}$ and the \emph{unit circle} $S:= \{x \in V : ||x|| = 1\}$ of $(V,\|\cdot\|)$. Thus $B$ is a compact, convex set centered at $o$, which is an interior point of it.
A normed plane is said to be \textit{strictly convex} if the triangle inequality is strict for vectors in distinct directions. We deal with three orthogonality types (see also the expository papers
\cite{Al-Be1}, \cite{Al-Be2}, and \cite{alonso}).
Two non-zero vectors $x,y \in (V,||\cdot||)$ are said to be\\

$\bullet$ \textit{Birkhoff orthogonal} whenever $||x+ty|| \geq ||x||$ for every $t \in \mathbb{R}$,  and in this case we write $x \dashv_B y$, \\

$\bullet$ \textit{Roberts orthogonal} if $||x+ty|| = ||x-ty||$ for all $t \in \mathbb{R}$,  denoted by $x \dashv_R y$, and\\

$\bullet$ \textit{isosceles orthogonal} when $||x+y|| = ||x-y||$, denoted by $x \dashv_I y$.\\

  It is worth mentioning that (uncommonly, but useful) we prefer to restrict our orthogonality definitions to non-zero vectors.
  For $x,y \in V$ we denote by $[xy]$, $\left<xy\right>$ and $\left.[xy\right>$ the \emph{closed line segment} connecting $x$ and $y$ (an \emph{open segment} is denoted by $(xy)$), the \emph{line} spanned by $x$ and $y$, and the \emph{half-line} with origin $x$ and passing through $y$, respectively. Given $p \in V$ and a line $r \subseteq V$, we denote the usual distance from $p$ to $r$ by $d(p,r):= \inf_{q\in r}||p-q||$.

  The final part of this introductory section consists of some needed elementary results from the geometry of Minkowski planes, all of them taken from \cite{martini1}; this means that for
 proofs and more details the reader is referred to \cite{martini1}. The first ones characterize strictly convex normed planes as the ones whose unit circle does not contain a nondegenerate line segment.

\begin{prop}\label{prop} Let $a, b, c \in V$ be three non-collinear points. Then, we have the equality $||a-c|| = ||a-b|| + ||b-c||$ if and only if there exists a segment $L\subseteq S$  containing the unit vectors $\frac{b-a}{||b-a||}$ and $\frac{c-b}{||c-b||}$. In this case, we also have $\frac{c-a}{||c-a||} \in \mathrm{relint}L$.
\end{prop}

\begin{lemma}\label{lemma2} Let $(V,||\cdot||)$ be a normed plane, and $y,z \in V$ be distinct points. Assume that $w \in (yz)$.
Then for every $x \in V$ we have $||x-w|| \leq \max(||x-y||,||x-z||)$, with equality if and only if $||x-w|| = ||x-y|| = ||x-z||$ . If equality holds, we have the following consequences: \\

\noindent\textbf{(a)} $||x-w||$ is the shortest distance from $x$ to the line $\left<yz\right>$.  \\

\noindent\textbf{(b)} The segment $[yz]$ is contained in the circle with center $x$ and radius $||x-w||$. Hence, $||x-v|| = ||x-w||$ for every $v \in [yz]$.  \\

\noindent In particular, the equality case cannot occur in strictly convex Minkowski planes.
\end{lemma}


The next lemma concerns distances from points to lines. This will be important in the study of bisectors of two points, and the last proposition is an easy consequence of the triangle inequality for quadrilaterals.

\begin{lemma} \label{lemma8} Let $r$ be a line in a normed plane and consider points $p \notin r$ and $q \in r$ such that $d(p,r) = ||p-q||$. If $p' \notin r$ and $q' \in r$ are points for which $\left<p'q'\right>$ is parallel to $\left<pq\right>$, then $d(p',r) = ||p'-q'||$.
\end{lemma}


\begin{prop}\label{prop4} Let \textbf{abcd} be a convex quadrilateral in a normed plane $(V,||\cdot||)$, with vertices in this written order. Then \\
\[ ||a-c|| + ||b-d|| \geq ||a-b|| + ||c-d||, \]\\
\noindent
with equality if and only if $[vw] \subseteq S$, where $v = \frac{c-a}{||c-a||}$, $w = \frac{b-d}{||b-d||}$, and $S$ is the unit circle. The same holds for the other pair of opposite sides. Notice that, in particular,
the sum of lengths of the diagonals cannot be equal to the sum of lengths of two opposite sides in a stricly normed plane.
\end{prop}


\section{The geometric structure of bisectors}\label{bisectors}

Given two distinct points $x,y \in (V,||\cdot||)$, we define the \textit{bisector} of $x$ and $y$ (or of the segment $[xy]$) to be the set \\
\[
\mathrm{bis}(x,y):= \{z \in V: ||z-x|| = ||z-y||\}.
\]\\
 Geometric properties of bisectors in arbitrary Minkowski planes can be, as is well known, quite complicated (see the surveys \cite{martini1} and \cite{martini2}; for bisectors in higher dimensional normed spaces we add the references \cite{horvath1} and \cite{horvath2}). In \cite{jahn} a general geometric description of bisectors is given. We say that a pair $(x,y)$ of vectors is a \textit{strict pair} if $[xy]$ is not parallel to a segment of the unit circle. Otherwise we say that $(x,y)$ is a \textit{non-strict pair}. Also, we define a \textit{cone} to be the convex hull of two half-lines with the same origin (called the \textit{apex} of the cone). In \cite{jahn} it is
proved that if $(x,y)$ is a strict pair, then $\mathrm{bis}(x,y)$ is a curve which is homeomorphic to a line, and if $(x,y)$ is a non-strict pair, then $\mathrm{bis}(x,y)$ is the union of two cones with a curve connecting
  their apices and itself homeomorphic to a closed interval. The present section is devoted to tackle this theory from another point of view (a little more geometric). This point of view will also be useful for
  characterizing Roberts orthogonality in the next section.

\begin{prop}\label{prop23} Let $(V,||\cdot||)$ be a Minkowski plane and let $x,y \in V$ be distinct points. Then any line $l$ parallel to $y - x$ intersects $\mathrm{bis}(x,y)$. This intersection is given by only one point
for any line $l$ if and only if (x,y) is a strict pair. In particular, a Minkowski plane is strictly convex if and only if for every pair of points $x,y \in V$ any line parallel to $y-x$ intersects $\mathrm{bis}(x,y)$ in exactly one point.
\end{prop}
\noindent
\textbf{Proof.} First, fix distinct points $x, y \in V$ and a line $l$ which is parallel to $y - x$. If $l = \left<xy\right>$, then the midpoint of $[xy]$ belongs to $l\cap \mathrm{bis}(x,y)$. Assume now that $l$ is not the line $\left<xy\right>$. Let $p_1 \in l$ be a point such that $||x-p_1|| = d(x,l)$. By Lemma \ref{lemma8} and taking $p_2 \in l$ such that $[yp_2]$ is parallel to $[xp_1]$, we have $d(y,l) = ||y-p_2| = ||x-p_1|| = d(x,l)$. Consider the continuous function $p \mapsto ||y-p|| - ||x-p||$ with $p$ ranging over $l$. This function is non-negative at $p_1$ and non-positive at $p = p_2$. Hence, by the Intermediate Value Theorem we have $||x-p|| = ||y-p||$ for some $p \in [p_1p_2] \subseteq l$. \\

\begin{figure}
\centering
\includegraphics{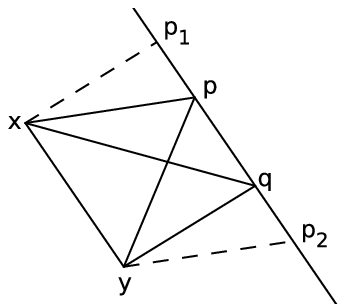}
\caption{Proposition \ref{prop23}}
\label{fig39}
\end{figure}

Assume now that there exists a line $l$ containing two distinct points $p$ and $q$, say, of $\mathrm{bis}(x,y)$. Hence, the points $x$, $y$, $p$ and $q$ are vertices of a quadrilateral for which the sum of lengths of the diagonals equals the sum of the lengths of two opposite sides. By Proposition \ref{prop4} it follows that (interchanging $p$ by $q$, if necessary) the segment $\left[\frac{x-q}{||x-q||} \frac{y-p}{||y-p||}\right]$ is contained in the unit circle. To show that this segment is parallel to $y - x$ it is enough to prove that $||x-q|| = ||y-p||$ (since $p-q$ and $y-x$ are parallel). If the segments $[xp]$ and $[yq]$ are parallel, this is obvious. Thus we suppose that they are not parallel. Choose points $p_1$ and $p_2$ in $l$ such that $[xp_1]$ is parallel to $[yq]$, and $[yp_2]$ is parallel to $[xp]$ (see Figure \ref{fig39}). By Lemma \ref{lemma2} it follows that $||x-p|| \leq \mathrm{max}(||x-p_1||,||x-q||) = ||y-q||$ and $||y-q|| \leq \mathrm{max}(||y-p_2||,||y-p||) = ||x-p||$. This shows the desired equality.\\

Suppose now that the unit circle contains a segment $[ab]$ which is parallel to $[xy]$. The line passing through the origin $o$ and the point $\frac{a-b}{2}$ is parallel to $b - a$, and it is easy to see that these two points belong to $\mathrm{bis}\left(a,\frac{a+b}{2}\right)$. Therefore, since the segment $\left[a\frac{a+b}{2}\right]$ is still parallel to $y - x$, it follows by translation and homothety that there exists a line $l$ parallel to $y - x$ containing more than one point of $\mathrm{bis}(x,y)$.

\begin{flushright} $\square$ \end{flushright}

\begin{remark}\normalfont The characterization of stricly convex norms given in the last proposition appeared for the first time in \cite{holub}.
\end{remark}

We will prove now that the bisector of a strict pair must be homeomorphic (in the induced topology of $V$) to a line. First we need an auxiliary lemma.

\begin{lemma}\label{lemma16} A segment $[xy] \subseteq V$ is not parallel to any segment of the unit circle if, and only if, for every  $z \in \mathrm{bis}(x,y) \setminus [xy]$ it holds that $\mathrm{bis}(x,y)$ is contained in the union of the convex region $\mathrm{conv}(\left.[zx\right>\cup\left.[zy\right>)$ and its image symmetric through $z$.
\end{lemma}
\noindent
\textbf{Proof.} If $[xy]$ is a segment which is parallel to some segment of the unit circle, then by the previous proposition there exists a line which is parallel to $\left<xy\right>$ containing (at least) two points of $\mathrm{bis}(x,y)$. Choosing one of them to be $z$, it is clear that the other one does not belong to the described region. \\

\begin{figure}
\centering
\includegraphics{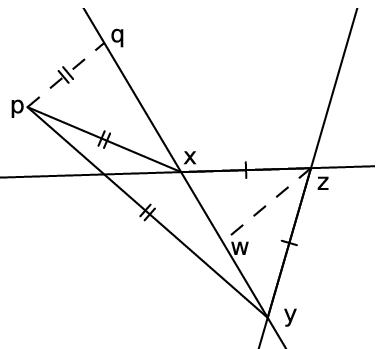}
\caption{Lemma \ref{lemma16}}
\label{fig40}
\end{figure}

Let now $[xy]$ be a segment which is not parallel to a segment of the unit circle, and $z \in \mathrm{bis}(x,y)\setminus [xy]$. Assume that there exists a point $p \in \mathrm{bis}(x,y)$ that does not lie in the described region. We have to consider two cases. First, suppose that the points $x,y,z$ and $p$ form a convex quadrilateral. Hence, we just have to draw the segments from $z$ to $\left<xy\right>$ which are respectively parallel to $[px]$ and $[py]$, and notice that both have the same length as $[zx]$ and $[zy]$ (see Lemma \ref{lemma2}). This contradicts the hypothesis on $[xy]$. If $x,y,z$ and $p$ do not form a convex quadrilateral, then we proceed as follows: by Lemma \ref{lemma2}, the distance from $z$ to $\left<xy\right>$ is particularly satisfied by some $w \in (xy)$. Also, by Lemma \ref{lemma8} the distance from $p$ to $\left<xy\right>$ must be attained, in particular, for some $q \in \left<xy\right>$ for which the segments $[zw]$ and $[pq]$ are parallel. Since $q \notin [xy]$, it follows that $||p-q|| = ||p-x|| = ||p-y||$, and hence Lemma \ref{lemma2} guarantees that $[xy]$ is parallel to a segment of the unit circle (see Figure \ref{fig40}).

\begin{flushright} $\square$ \end{flushright}

\begin{prop}\label{prop19} If $[xy] \subseteq V$ is a segment which is not parallel to any segment of the unit circle, then $\mathrm{bis}(x,y)$ is homeomorphic (in the topology induced by $V$) to a line.
\end{prop}
\noindent
\textbf{Proof.} Since $\mathrm{bis}(x,y)$ is symmetric with respect to $\frac{x+y}{2}$, it is clearly sufficient to prove that the intersection of $\mathrm{bis}(x,y)$ with one of the half-planes determined by $\left<xy\right>$ is homeomorphic to the interval $[0,\infty)$. In view of Proposition \ref{prop23} we may define a function $p:[0,\infty) \rightarrow V$ that associates each non-negative number $d$ to the point $p(d) \in \mathrm{bis}(x,y) \cap l(d)$, where $l(d)$ is the line parallel to $\left<xy\right>$ at distance $d$. We shall show that $p$ is a homeomorphism over its image. Notice that $p$ is injective (Proposition \ref{prop23}), and it is also clear from the continuity of the distance function that its inverse is continuous. Thus, we just have to show that $p$ is continuous. We start by showing that $p$ is continuous at any $d_0 > 0$. Fix such a number and let $(d_n)_{n \in \mathbb{N}}$ be a sequence converging to $d_0$. If $n_0$ is a natural number such that $d_n \leq d_0+1$ whenever $n > n_0$, then by Lemma \ref{lemma16} we have that the set $\{p(d_n)\}_{n > n_0}$ is contained in the compact set $\mathrm{conv}\{x,y,p(d_0)\}\cup\mathrm{conv}\{x_1,y_1,p(d_0)\}$, where $x_1$ and $y_1$ are the intersections of $l(d_0+1)$ with the lines $\left<xp(d_0)\right>$ and $\left<yp(d_0)\right>$, respectively. Hence the sequence $(p(d_n))_{n\in \mathbb{N}}$ is contained in some compact set of $V$, and then it has a converging subsequence $(p(d_{n_k}))_{k\in \mathbb{N}}$. Since $d_{n_k} \rightarrow d_0$, it follows that $p(d_{n_k})$ converges to a point of $l(d_0)$. Moreover, since $||x-p(d_{n_k})|| = ||y-p(d_{n_k})||$ for every $k \in \mathbb{N}$, we have that the limit point of $p(d_{n_k})$ belongs to $\mathrm{bis}(x,y)$. Therefore $p(d_{n_k}) \rightarrow p(d_0)$. Notice that the same argument shows that any converging subsequence of $(p_n)_{n\in\mathbb{N}}$ must converge to $p(d_0)$. Thus, by standard analysis it follows that $p(d_n)$ converges itself to $p(d_0)$. To prove that $p$ is continuous at $d = 0$, we repeat the argument and see that if $d_n \rightarrow 0$, then for some $n_0 \in \mathbb{N}$ it holds that $p(d_n) \in \mathrm{conv}\{x,y,p(1)\}$ whenever $n > n_0$.  \\

\begin{flushright} $\square$ \end{flushright}

Now we study the geometric structure of bisectors for non-strict pairs. This is established in the next two propositions. But first we notice that if $(x,y)$ is such a pair, then the segment $[xy]$, is a maximal segment of precisely two circles of the plane, each of them with its center lying in one of the half-planes determined by the line $\left<xy\right>$. In fact, if $[xy]$ is parallel to a maximal segment $[ab] \subseteq S$ and $\lambda = \frac{||x-y||}{||a-b||}$, then (assuming, without loss of generality, that $y - x$ = $\lambda b - \lambda a$) we just have to consider the circles $\lambda S + x - \lambda a$ and $\lambda S + x + \lambda b$.

\begin{figure}
\centering
\includegraphics{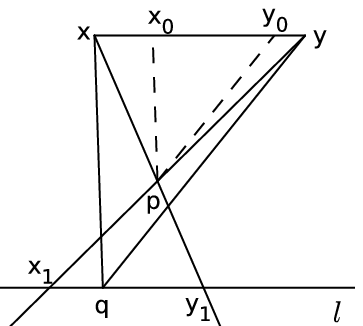}
\caption{$l \cap \mathrm{bis}(x,y) = [x_1y_1]$}
\label{fig41}
\end{figure}

\begin{prop}\label{prop21} Let $[xy] \subseteq V$ be a segment which is parallel to a segment of the unit circle, and let $p \in V$ be the center of one of the (two) circles which contain $[xy]$ as maximal segment. Let $l$ be any line parallel to $\left<xy\right>$ such that $p$ lies in the interior of the strip determined by these two lines. Then $l \cap \mathrm{bis}(x,y)$ is precisely the segment $[x_1y_1]$, where $x_1$ and $y_1$ are the intersections of $\left.[xp\right>$ and $\left.[yp\right>$ with $l$, respectively. In particular, $\mathrm{bis}(x,y)$ contains the cone $\mathrm{conv}(\left.[p(2p-x)\right>\cup\left.[p(2p-y)\right>)$ and, consequently, $\mathrm{int}(\mathrm{bis}(x,y)) \neq \emptyset$.
\end{prop}
\noindent
\textbf{Proof.} Let $q \in [x_1y_1]$. It is immediate that the parallels to $\left<qx\right>$ and $\left<qy\right>$ through $p$ intersect the line $\left<xy\right>$ in the interior of the segment $[xy]$ (see Figure \ref{fig41}). Let $x_0$ and $y_0$ be these intersection points, respectively. Since $[xy]$ belongs to a circle with center $p$, it follows that $||p-x_0|| = ||p-y_0||$. Also, it is clear that $\frac{||q-x||}{||p-x_0||} = \frac{||q-y||}{||p-y_0||}$ (the endpoints of the segments) correspondingly lie in the same parallel lines, and then we have $||q-x|| = ||q-y||$. \\

Assume now that $q \in \left<x_1y_1\right>\setminus [x_1y_1]$. We may suppose that, without loss of generality, $y_1$ lies between $q$ and $x_1$. If $x_0$ and $y_0$ are chosen in $\left<xy\right>$ such that $[px_0]$ is parallel to $[qx]$ and $[py_0]$ is parallel to $[qy]$, then it is easy to prove that $x_0 \in [xy]$ but $y_0 \notin [xy]$. Since $[xy]$ is maximal, we therefore must have $||p-y_0|| > ||p-x_0||$ (Lemma \ref{lemma2}).
Hence $||q-y|| > ||q-x||$.
\begin{flushright} $\square$ \end{flushright}

We described the shape of bisector $\mathrm{bis}(x,y)$ outside the strip given by the lines parallel to $\left<xy\right>$ and passing through the centers of the circles which contain $[xy]$ as maximal segment.
We now describe the shape of $\mathrm{bis}(x,y)$ within this strip.

\begin{figure}
\centering
\includegraphics{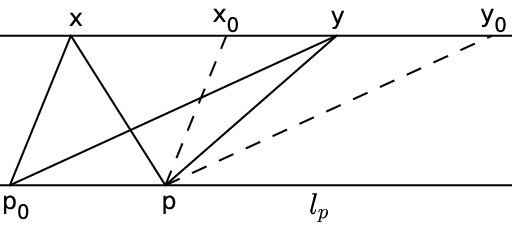}
\caption{$l=l_p$}
\label{fig46}
\end{figure}

\begin{prop}\label{prop22} Let $[xy]$ be a segment which is parallel to a segment from the unit circle, and assume that $p$ and $q$ are the centers of the circles which contain $[xy]$ as a maximal segment. Let $l_p$ and $l_q$ be the parallels to $\left<xy\right>$ through $p$ and $q$, respectively. Then the intersection $\mathrm{bis}(x,y)\cap\mathrm{conv}(l_p\cup l_q)$ is a curve from $p$ to $q$ which is homeomorphic to a compact interval.
\end{prop}
\noindent\textbf{Proof.} Notice that since $\mathrm{bis}(x,y)$ is symmetric with respect to $\frac{x+y}{2}$, it is enough to prove the result for the strip determined by $\left<xy\right>$ and $l_p$. First we show that if $l \subseteq \mathrm{conv}(\left<xy\right>\cup l_p)$ is a line parallel to $\left<xy\right>$, then the intersection $l\cap \mathrm{bis}(x,y)$ contains precisely one point. This is obvious if $l = \left<xy\right>$.
Assume now that $l = l_p$ and suppose that there exists a point $p_0 \neq p$ with $p_0 \in l_p\cap \mathrm{bis}(x,y)$. Let $x_0,y_0 \in \left<xy\right>$ be such that $[px_0]$ is parallel to $[p_0x]$, and $[py_0]$ is parallel to $[p_0y]$. Then, renaming the points if necessary, we may say that $x_0 \in (xy)$ and $y_0 \in \left<xy\right>\setminus [xy]$ (see Figure \ref{fig46}). Since $[xy]$ belongs to a circle with center $p$, we have that $||p-x_0|| = ||p-x||$. In particular, it follows that $||p-x|| = ||p-y|| = ||p-y_0||$. Thus, by Lemma \ref{lemma2} we see that the segment $[xy_0]$ belongs to a circle with center $p$, and this contradicts the maximality of $[xy]$. \\

\begin{figure}
\centering
\includegraphics{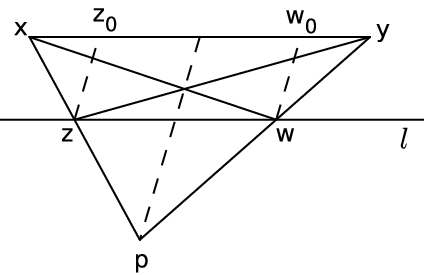}
\caption{$l \subseteq \mathrm{int}(\mathrm{conv}(\left<xy\right>\cup l_p))$}
\label{fig47}
\end{figure}

Let now $l$ be a line parallel to and strictly between $\left<xy\right>$ and $l_p$. Let $z$ and $w$ be the intersections of $l$ with the segments $[xp]$ and $[yp]$, respectively. It is easy to see that the segments $[zx]$, $[wy]$, $[zz_0]$, and $[ww_0]$ have the same length, where $z_0$ and $w_0$ are the points of $[xy]$ such that $[zz_0]$ and $[ww_0]$ are parallel to $\left[p\frac{x+y}{2}\right]$ (see Figure \ref{fig47}). By Lemma \ref{lemma2} it follows that $||z-y|| > ||z-x||$ and $||w-x|| > ||w-y||$; one may wonder whether equality cannot hold. If this would happen, then $[xy]$ would be a segment of a circle with radius $||p-z|| < ||p-x||$, and this contradicts the maximality of $[xy]$ in the circle with radius $||p-x||$. Now, by the Intermediate Value Theorem applied to the function $v \in [zw] \mapsto ||v-x|| - ||v-y||$, it follows that there exists a point $v_0 \in (zw)\cap \mathrm{bis}(x,y)$. To prove that $v_0$ is the only point of $l\cap \mathrm{bis}(x,y)$, we just have to repeat the proof for the case $l = l_p$ (actually, we would have a segment properly containing $[xy]$ in a circle with radius smaller than $||p-x||$). \\

From the previous argument it also follows that $\mathrm{bis}(x,y) \cap \mathrm{conv}(\left<xy\right>\cap l_p)$ is contained in the compact set $\mathrm{conv}\{x,y,p\}$. Hence, we may repeat the proof of Proposition \ref{prop19} to prove that the function which associates each $d \in [0,||p-x||]$ to the point $l_d \cap \mathrm{bis}(x,y)$, where $l_d$ is the parallel to $\left<xy\right>$ at distance $d$, is a homeomorphism.
\begin{flushright} $\square$ \end{flushright}

\begin{coro}\label{coro14} The bisector of a segment $[xy]$ can contain at most one line. If $[xy]$ is parallel to some segment of the unit circle, then this line must be necessarily the line through the centers of the two circles which contain $[xy]$ as maximal segment.
\end{coro}
\noindent\textbf{Proof.} If $(x,y)$ is a strict pair, then the assertion follows immediately from Proposition \ref{prop19}. For the other case, notice that if $\mathrm{bis}(x,y)$ contains a line, then the curve described in Proposition \ref{prop22} must be a segment passing through the mentioned points.

\begin{flushright} $\square$ \end{flushright}

\begin{remark}\label{remark7} \normalfont By homothety and translation it follows that if the bisector $\mathrm{bis}(x,y)$ of a segment $[xy]$ contains a line $l$, then the bisector of any parallel segment also contains a line (in the same direction as $l$).
\end{remark}

\begin{figure}[h]
\centering
\includegraphics{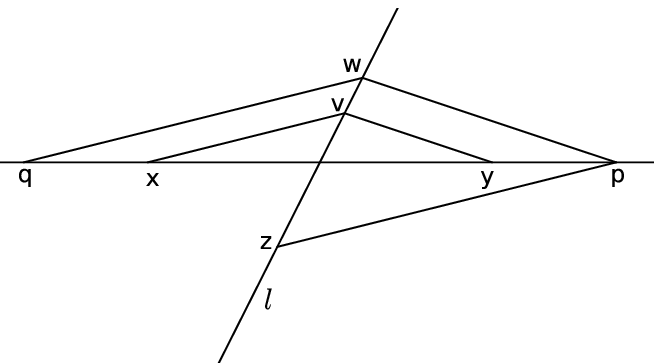}
\caption{Lemma \ref{lemma18}}
\label{fig50}
\end{figure}

\begin{lemma}\label{lemma18} Let $[xy]\subseteq V$ be a segment and assume that there exists a line $l$ contained in $\mathrm{bis}(x,y)$. Then, for any segment $[zw] \subseteq l$ centered at $\frac{x+y}{2}$, we have $\left<xy\right> \subseteq \mathrm{bis}(z,w)$.
\end{lemma}
\noindent\textbf{Proof.} Let $[zw]$ be such a segment and fix an arbitrary point $p \in \left<xy\right>$ which lies in the same half-plane determined by $l$ as $y$. Assume that $q \in \left<xy\right>$ is symmetric to $p$ through $\frac{x+y}{2}$. Hence $||z-p|| = ||w-q||$. Let $v \in l$ be such that $[vy]$ and $[wp]$ are parallel segments and notice that, thus, $[xv]$ is parallel to $[wq]$ (see Figure \ref{fig50}). Since $l \subseteq \mathrm{bis}(x,y)$, it follows that $||v-x|| = ||v-y||$. Then, by homothety, we have $||w-q|| = ||w-p||$. Therefore $||z-p|| = ||w-p||$, as we wished.

\begin{flushright} $\square$ \end{flushright}

\section{New characterizations of Roberts orthogonality}\label{roberts}

\noindent Considering homothety and translation, one can observe that the geometric structure of bisectors in normed planes can be studied by looking at the bisectors of the diameters of the unit circle. Despite the simplicity of the proof, the next theorem is important. Namely, it justifies why the geometric constants, which will be defined later, indeed quantify the difference between Roberts orthogonality and Birkhoff orthogonality.

\begin{teo}\label{teoroberts} Let $x,y \in S$ be linearly independent unit vectors. Then the following statements are equivalent:\normalfont \\

\noindent\textbf{(a)} \textit{The bisector $\mathrm{bis}(-x,x)$ contains the line $\left<oy\right>$}.\\

\noindent\textbf{(b)} \textit{$x \dashv_R y$}. \\

\noindent\textbf{(c)} \textit{The unit circle is invariant through the linear transformation $T:V \rightarrow V$ defined by setting $T(x) = x$ and $T(y) = -y$.}
\end{teo}

\noindent
\textbf{Proof.} If \textbf{(a)} holds, then $||ty - x|| = ||ty + x||$ for every $t \in \mathbb{R}$, and this means that $x \dashv_R y$. If \textbf{(b)} holds, then we have that $||\alpha x + \beta y|| = ||\alpha x - \beta y||$ for every $\alpha,\beta \in \mathbb{R}$. Hence $T$ is an isometry, and therefore \textbf{(c)} follows. Assume now that \textbf{(c)} is true. Since $T$ is an isometry, we have that $||ty-x|| = ||T(ty-x)|| = ||ty+x||$ for every $t\in\mathbb{R}$, and this gives $\left<oy\right> \subseteq \mathrm{bis}(-x,x)$.

\begin{flushright} $\square$ \end{flushright}

Using this theorem to study bisectors of chords in the unit circle we are able to provide a coordinate-free characterization of ellipses among all centrally symmetric two-dimensional convex figures.

\begin{prop}\label{propchords} Let $[xy]$ be a chord of the unit circle. Then $\mathrm{bis}(x,y)$ contains a line if and only if the unit circle is invariant with respect to the linear transformation $T:V \rightarrow V$ such that $T(y-x) = y-x$ and $T(x+y) = -x-y$.
\end{prop}
\noindent\textbf{Proof.} If $\mathrm{bis}(x,y)$ contains a line, then it must be necessarily the line $l$ passing through the origin $o$ and the midpoint $w = \frac{x+y}{2}$ (see Lemma \ref{lemma18}). Let $[(-v)v]$ be the diameter of the unit circle which is parallel to $[xy]$. It is clear that $\mathrm{bis}(-v,v)$ must contain a line parallel to $l$. But since the origin is obviously contained in the bisector of $[(-v)v]$, we have $l \subseteq \mathrm{bis}(-v,v)$. Now, from Theorem \ref{teoroberts} it follows that the unit circle is invariant through the linear map $T:V\rightarrow V$ for which $T(v) = v$ and $T(w) = - w$, and this clearly yields the assertion. For the converse, assume again that $[(-v)v]$ is the diameter of the unit circle which is parallel to $[xy]$, and let $w = \frac{x+y}{2}$. Clearly, the hypothesis gives that the unit circle is invariant with respect to the linear transformation $T:V \rightarrow V$ for which $T(v) = v$ and $T(w) = -w$. Hence $\mathrm{bis}(-v,v)$ contains a line, and therefore also $\mathrm{bis}(x,y)$ does.

\begin{flushright} $\square$ \end{flushright}

\begin{figure}
\centering
\includegraphics{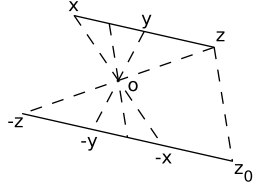}
\caption{$[(-z)(-x)] \subset [-zz_0]$}
\label{fig54}
\end{figure}

\begin{lemma} \label{lemmabis} If every bisector in a normed plane $(V,||\cdot||)$ contains a line, then any bisector in $V$ is, in fact, a line.
\end{lemma}
\noindent\textbf{Proof.} We just have to prove that if every bisector contains a line, then the plane is strictly convex. Assume the hypothesis and suppose that $[xz]$ is a segment of the unit circle. We can take $[xz]$ to be maximal. Let $y$ be the midpoint of $[xz]$. Since the bisector of $[xy]$ contains a line, it follows that the unit circle is invariant under the linear transformation $T$ which takes $T(x+y)$ to $-x-y$ and $T(y-x)$ to $y-x$. Then we have $z_0 = T(z) \in S$. Writing $z = 2y - x$, we have $z_0 = T(2y-x) = 2T(y)-T(x) = 2(-x) + (-1)(-y)$. It follows that $-x$ is the midpoint between $z_0$ and $(-y)$. In particular, the segment $[-zz_0]$ properly contains the segment $[(-z)(-x)]$  (see Figure \ref{fig54}). This contradicts the maximality of $[xz]$.

\begin{flushright} $\square$ \end{flushright}

\begin{coro}\label{coroellipse} A centrally symmetric convex body $K$ contained in a two-dimensional vector space is an ellipse if and only if for every $x,y \in \partial K$ it holds that $K$ is invariant under the linear transformation $T$ defined by setting $T(x+y) = -x-y$ and $T(y-x) = y-x$.
\end{coro}
\noindent
\textbf{Proof.} It is known that a normed plane is an inner product plane if and only if the bisector of each segment is a line (see \cite{martini2}). Thus, the assertion follows immediately from Proposition \ref{propchords} and Lemma \ref{lemmabis}.
\begin{flushright} $\square$ \end{flushright}

We finish this section with a characterization of inner product planes that we will need later, too.

\begin{prop}\label{propchar} A norm is derived from an inner product if and only if for every $x \in S$ the bisector $\mathrm{bis}(-x,x)$ contains the line segment connecting the points of $S$ at which the direction $x$ supports the unit ball.
\end{prop}
\noindent
\textbf{Proof.} Assume that $(V,||\cdot||)$ is not an inner product space. Then there exist $x,y\in S$ such that $y \dashv_B x$ but $||x-y|| \neq ||x+y||$ (see \cite{alonso}, Theorem 5.1). Thus, the segment joining the points of $S$ where $x$ supports the unit ball is $[(-y)y]$, but $y \notin \mathrm{bis}(-x,x)$. The converse is trivial.

\begin{flushright} $\square$ \end{flushright}

\section{The constant $c_B$}\label{cb}

Let $x \in S$ be a unit vector. We define the \textit{inner bisector} of the segment $[(-x)x]$ to be the set $\mathrm{bis_I}(-x,x) := \mathrm{bis}(-x,x)\cap B$. In other words, the inner bisector of a diameter is the set of points of its bisector which lie in the unit ball. We also define the \textit{inner projection} of $\mathrm{bis}(-x,x)$ to be the set \\
\[ \mathrm{P_I}(x) := \left\{\frac{z}{||z||}:z \in \left(\mathrm{bis}(-x,x)\cap B\right) \setminus \{o\}\right\}. \]\\
The intuitive reason why we define these sets is the following: fix $x \in S$ and assume that $y \in S$ is such that $y \dashv_B x$. If $\mathrm{bis}(-x,x)$ contains a line, then it is clear that $\mathrm{P_I}(x) = \{-y,y\}$ (recall that Roberts orthogonality implies Birkhoff orthogonality; cf. \cite{alonso}). If this is not the case, then $x$ is not Roberts orthogonal to $y$, and the inner projection cannot be a two-point set anymore. Measuring, somehow, how far the inner projection is from $y$ is a way to quantify the difference between Birkhoff orthogonality and Roberts orthogonality. This will be made more precise after a few more steps (using a generalized sine function). Related properties of radial projections of bisectors were studied in \cite{mw}.

\begin{prop}\label{propinnerbis} For any $x \in S$, the inner bisector $\mathrm{bis_I}(-x,x)$ is a curve which is homeomorphic to a closed interval.
\end{prop}
\noindent\textbf{Proof.} This follows immediately from Propositions \ref{prop19}, \ref{prop21}, and \ref{prop22}. Notice that if $(-x,x)$ is a non-strict pair, then the inner bisector is a portion of the 1-dimensional component of $\mathrm{bis}(-x,x)$, and it equals all of it if and only if the plane is rectilinear.

\begin{flushright} $\square$ \end{flushright}

\begin{figure}
\centering
\includegraphics{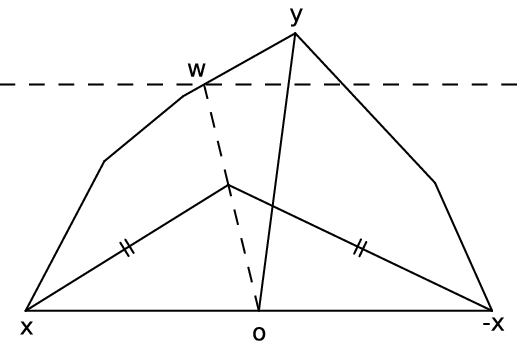}
\caption{$s(w,x), \ w\in \mathrm{P_I(x)}$}
\label{fig1e}
\end{figure}

In \cite{bmt} the function $s:S\times S \rightarrow \mathbb{R}$ given by $s(x,y) = \inf\{||x+ty|| : t \in \mathbb{R}\}$ is studied. This function somehow plays the role of the sine function in a normed plane (or space). It is known that $s(x,y) \leq 1$ for every $x,y \in S$, and that equality holds if and only if $x \dashv_B y$ (see again \cite{bmt}). Having said this, we define the constant $c_B(||\cdot||)$ to be \\
\[ c_B(||\cdot||) := \inf_{x\in S}\left(\inf_{w\in \mathrm{P_I}(x)}s(w,x)\right). \] \\
Geometrically we have that if $y \dashv_B x$, then the value of $s(w,x)$ is the length of the segment whose endpoints are the origin and the intersection between the line $t\mapsto w+tx$ and the segment $[oy]$. Hence, the number $\inf_{w\in \mathrm{P_I}(x)}s(w,x)$ is the infimum of the lengths that the parallels to $x$ passing through the points of $\mathrm{P_I(x)}$ determine over the segment $[oy]$. Figure \ref{fig1e} illustrates the situation. Before studying upper and lower bounds for $c_B$, we shall calculate this constant, as an example, for rectilinear planes.

\begin{example}\label{examplerect}\normalfont Choose $(V,||\cdot||)$ as a rectilinear plane, i.e., its unit circle be a parallelogram. Then $c_B(||\cdot||) = \frac{1}{2}$.
\end{example}
\noindent\textbf{Proof.} Let $(V,||\cdot||)$ be a rectilinear plane and assume that $p$ and $q$ are consecutive vertices of the unit circle. It is clear that we can consider $x$ ranging through the segment $\left[p\left(\frac{p+q}{2}\right)\right]$ to perform our calculations. The first step is to determine the structure of the inner bisectors. Notice that $\mathrm{bis_I}(-p,p) = [(-q)q]$
and bis$_I\!\left(-\frac{p+q}{2},\frac{p+q}{2}\right) = \left[-\left(\frac{q-p}{2}\right)\frac{q-p}{2}\right]$. It follows immediately that $\inf_{w\in\mathrm{P_I}(x)}s(w,x) = 1$ if $x = p$ or $x = \frac{p+q}{2}$. \\

\begin{figure}[h]
\centering
\includegraphics{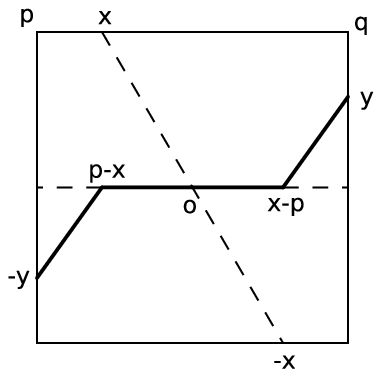}
\caption{$\mathrm{bis_I}(-x,x)$}
\label{fig2e}
\end{figure}

Now let $x \in \left(p\frac{p+q}{2}\right)$, and let $y \in [q(-p)]$ be such that $||y - q|| = ||x - p||$. Then it is easy to see that $\mathrm{bis_I}(-x,x)$ is the union of the segments $[(-y)(p-x)]$, $[(p-x)(x-p)]$ and $[(x-p)y]$ (see Figure \ref{fig2e}). Hence the inner projection $\mathrm{P_I}(x)$ is the union of the segment $\left[y\left(\frac{q-p}{2}\right)\right]$ with its symmetric image. It follows from the geometric approach given above that $\inf_{w\in\mathrm{P_I}(x)}s(w,x)$ is attained for $w = \frac{q-p}{2}$ whenever $x \in \left(p\frac{p+q}{2}\right)$. A simple calculation gives $s\!\left(\frac{q-p}{2},x\right) = \frac{1}{2-||x-p||}$. Therefore, $c_B(||\cdot||) = \frac{1}{2}$.

\begin{flushright} $\square$ \end{flushright}

An interesting phenomenon appears in the previous example: the segments which form $\mathrm{P_I}(x)$ when $x$ ranges within $\left(p\frac{p+q}{2}\right)$ degenerate to points when $x = p$ or $x = \frac{p+q}{2}$. For that reason, the transformation which maps each $x \in S$ to the number $\inf_{w\in\mathrm{P_I}(x)}s(w,x)$ is not necessarily continuous. Moreover, the rectilinear plane is an example of a normed plane for which $c_B(||\cdot||)$ is not attained for a pair $x \in S$, $w \in \mathrm{P_I}(x)$.\\

 It seems to be difficult to calculate the constant $c_B$ for more complicated norms without computational methods. Nevertheless, we can give an easy and sharp upper bound and a lower bound (which, possibly, is not sharp, as will become clear in the proof).

\begin{teo}\label{teocb} Let $(V,||\cdot||)$ be a normed plane. Then \\
\[ \frac{1}{3} \leq c_B(||\cdot||) \leq 1,\] \\
and equality on the right holds if and only if the norm is derived from an inner product.
\end{teo}
\noindent\textbf{Proof.} The inequality $c_B(||\cdot||) \leq 1$ is obvious since $s(x,y) \leq 1$ for any $x,y \in S$. If equality holds, then for every $x \in S$ and $w \in \mathrm{P_I}(x)$ we have $s(w,x) = 1$. It follows that for every $x \in S$ the set $\mathrm{bis_I}(-x,x)$ is the segment connecting the two points of $S$ where $x$ supports $B$. From Proposition \ref{propchar} we have that $(V,||\cdot||)$ is an inner product plane. \\

To show the other inequality, we need first an auxiliary result: let $x \in S$ be an arbitrary unit vector and assume that $y \in S$ is such that $y \dashv_B x$. Denote by $H$ the (closed) half-plane determined by the line $\left<(-x)x\right>$ which contains $y$. Then $H\cap\mathrm{bis_I}(-x,x) \subseteq \mathrm{conv}\{o,y+2x,y-2x\}$. If the opposite holds, then we may assume, without loss of generality, that $H\cap\mathrm{bis_I}(-x,x)$ contains a point $p$ which lies in $\mathrm{int}\left(\mathrm{conv}\{o,x,y+2x\}\right)$ (indeed, by convexity and Birkhoff orthogonality we have that the lines $\left<y(y+2x)\right>$ and $\left<x(y+2x)\right>$ support $H\cap B$). We may write $p = \alpha x + \beta(y+2x)$ for some $\alpha,\beta > 0$ with $\alpha + \beta < 1$. Then the triangle inequality gives \\
\[ ||p-x|| = ||(\alpha+2\beta-1)x + \beta y|| \leq \beta + |\alpha + 2\beta - 1| \ \mathrm{and} \]

\[ ||p+x||  = ||(\alpha+2\beta+1)x+\beta y|| \geq \alpha + \beta + 1. \]\\
Thus, since $||p-x|| = ||p+x||$, it follows that $\alpha +1 \leq |\alpha + 2\beta - 1|$, and this implies $\beta \geq 1$ or $\alpha + \beta \leq 0$. In both cases we have a contradiction.\\

\begin{figure}[t]
\centering
\includegraphics{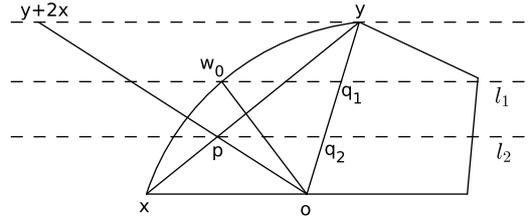}
\caption{$||q_1|| \geq ||q_2||$}
\label{fig4e}
\end{figure}

Now, assume that $w \in \mathrm{bis_I}(-x,x)\setminus\{o\}$, and let $l_1$ be the line parallel to $x$ and passing through $w_0 = \frac{w}{||w||}$. Let $p$ be the intersection of the segments $[o(y+2x)]$ and $[xy]$, and let $l_2$ be the line parallel to $x$ drawn through $p$. It is clear that $l_2$ intersects $[oy]$ in a point $q_2$ closer to the origin than the intersection point $q_1$  of $l_1$ with the same segment (see Figure \ref{fig4e}). It is easy to see that $||q_2|| = \frac{1}{3}$, and then the desired lower bound follows.

\begin{flushright} $\square$ \end{flushright}

\begin{remark}\normalfont If $(-x,x)$ is a strict pair, then the location of the inner bisector proved above is an easy consequence of Theorem 2.4 in \cite{vaisala}.
\end{remark}

\begin{open}\normalfont What is the sharp lower bound for $c_B(||\cdot||)$? Is it attained by some Minkowski plane?
\end{open}

\begin{prop}\label{propest} In any normed plane $(V,||\cdot||)$ we have the inequality \\
\[ c_B(||\cdot||) \leq D(||\cdot||), \]\\
where $D$ is the constant defined in \cite{wu} as $D:=\inf\left\{\inf_{t\in\mathbb{R}}||x+ty||:x,y\in S \ \mathrm{and} \ x\dashv_I y\right\}$.
\end{prop}
\noindent\textbf{Proof.} It is clear that we may write $D = \inf\{s(x,y):x,y \in S \ \mathrm{and} \ x\dashv_I y\}$. If $x \dashv_I y$, then $y \in \mathrm{P_I}(x)$. Hence $\inf_{w\in\mathrm{P_I}(x)}s(w,x) \leq s(y,x)$, and the desired follows.

\begin{flushright} $\square$ \end{flushright}

To finish this section, we use the same method as in Example \ref{examplerect} to calculate the constant $c_B$ for regular $(4n)$-gonal norms. These are the norms whose unit ball is an affine regular $(4n)$-gon.

\begin{figure}[h]
\centering
\includegraphics{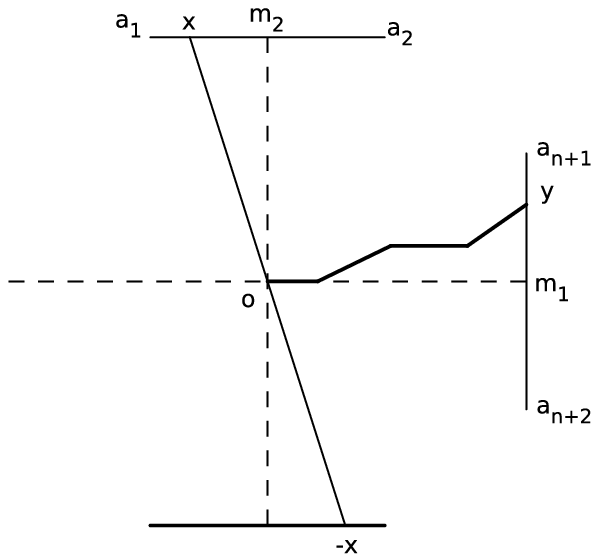}
\caption{$\mathrm{bis_I}(-x,x)$}
\label{fig3e}
\end{figure}

\begin{prop} Given $n \in \mathbb{N}$, let $(V,||\cdot||_{4n})$ denote a Minkowski plane whose unit circle is an affine regular $(4n)$-gon. Then \\
\[c_B\left(||\cdot||_{4n}\right) = \left(\cos\frac{\pi}{4n}\right)^2\,. \]
\end{prop}
\noindent\textbf{Proof.} The reason why it is not difficult to calculate the constant $c_B$ for regular $(4n)$-gonal norms is that we can locate the inner projections of its bisectors. Let $\mathbf{a_1a_2...a_{4n}}$ be an affine regular $(4n)$-gon which is the unit circle of the Minkowski plane $(V,||\cdot||_{4n})$ and denote by $m_2$ and $m_1$ respectively the midpoints of the sides $[a_1a_2]$ and $[a_{n+1}a_{n+2}]$ respectively. Then, $V$ is a \textit{symmetric Minkowski plane} and $\{m_1,m_2\}$ is a pair of axes, i.e., \\
\[ ||m_1 + tm_2|| = ||m_1 - tm_2|| = ||m_2+tm_1|| = ||m_2 - tm_1|| \]\\
for every $t \in \mathbb{R}$ (for more on symmetric Minkowski planes we refer the reader to \cite{wu}). It is clear that we may consider $x$ ranging through the segment $[a_1m_2]$ to describe all inner projections. It is also clear that if $x = a_1$ or $x = m_2$, then $\mathrm{bis_I}(-x,x)$ is a straight segment, and hence $\inf_{w\in\mathrm{P_I}(x)}s(w,x)=1$. We describe now the inner projection in the case $x \in (a_1m_2)$. We will consider only one of the half-planes determined by $\left<(-x)x\right>$ (namely, the one containing $m_1$, which we call $H$), since the bisector is symmetric through the origin. From \cite{ma}, Section 2.2, we have that $H\cap\mathrm{bis_I}(-x,x)$ must be a polygonal chain and it is easy to see that $[o(x-a_1)]$ is its first segment. Moreover, following \cite{wu}, Theorem 10, it is immediate that the point $y \in [a_{n+1}a_{n+2}]$ such that $||y - a_{n+1}|| = ||x-a_1||$ belongs to $\mathrm{bis_I}(-x,x)$. Since our polygonal circle is regular, we have that there exist precisely two directions in which the segments of $\mathrm{bis_I}(-x,x)$ can lie. One of them is the direction $m_1$ and the other one is the direction $a_{n+1}$ (see Figure \ref{fig3e}). It follows that the inner projection $\mathrm{P_I}(x)$ is precisely the segment $[m_1y]$ (if $\mathrm{bis_I}(-x,x)$ cuts the unit circle in a segment which is in the direction $a_{n+1}$) or it is contained in the segment $[m_1a_{n+1}]$ (otherwise), and thus $\inf_{w\in\mathrm{P_I}(x)}s(w,x)$ is attained for $w = m_1$. Indeed, since the direction $x$ supports the polygon at the vertex $a_{n+1}$ it follows that we can compute $s(w,x)$ looking at the distance from the intersection of the line $t \mapsto w + tx$ with the segment $[oa_{n+1}]$ to the origin. It is easy to see that this distance increases as $w$ ranges from $m_1$ to $y$.

Summarizing, if $x \in (a_1m_2)$, then $\inf_{w\in\mathrm{P_I}(x)}s(w,x) = s(m_1,x)$. Hence, in order to determine the constant $c_B$ we must calculate the infimum of the values of $s(m_1,x)$ as $x$ ranges within the segment $(a_1m_2)$. It is clear that the function $x \mapsto s(m_1,x)$ is increasing as $x$ goes from $a_1$ to $m_2$ (one can check this by using the same geometric argument used right above). Finally, since the sine function is continuous (see \cite{bmt}), we have\\
\[c_B\left(||\cdot||_{4n}\right) = \lim_{x\rightarrow a_1}s(m_1,x) = s(m_1,a_1). \]

To calculate $s(m_1,a_1)$ we take, for simplicity, the unit circle as the standard regular $(4n)$-gon in the Euclidean plane. The angle between the segments $[oa_{n+1}]$ and $[om_1]$ is $\frac{\pi}{4n}$, and the line parallel to $a_1$ drawn through $m_1$ cuts perpendicularly the segment $[oa_{n+1}]$ in a point $q$, say. Therefore, the value of $s(m_1,a_1)$ is the ratio between the Euclidean lengths of the segments $[oq]$ and $[oa_{n+1}]$. Basic trigonometry gives the desired value.

\begin{flushright} $\square$ \end{flushright}

\begin{remark}\normalfont One may wonder if it would be better to define $c_B(||\cdot||)$ to be $1 - \inf_{x\in S}\left(\inf_{w\in\mathrm{P_I}(x)}s(w,x)\right)$, replacing our definition in this way. This aesthetic change would not make any difference in the theory presented here, but it is worth mentioning that this constant would coincide with the constant $c_R$ defined in \cite{bmt} (and used to estimate how far a normed plane is from being Radon) for any regular $(4n)$-gonal norm. In some sense, this means that the difference between Roberts orthogonality and Birkhoff orthogonality in such a plane is as large as the plane is far from being Radon.
\end{remark}

\section{The constant $c_S$} \label{cs}

As we saw in Theorem \ref{teoroberts}, Roberts orthogonality is related to the invariance of the unit circle with respect to certain linear reflections (i.e., automorphisms of $V$ whose eigenvalues are $1$ and $-1$). Hence we can quantify the difference between Roberts orthogonality and Birkhoff orthogonality by estimating the distortion of the images of the unit circle with respect to the linear reflections whose eigenvectors are Birkhoff orthogonal vectors. Given two linearly independent vectors $x,y \in V$, denote by $T_{xy}$ the linear transformation defined by setting $T_{xy}(x) = x$ and $T_{xy}(y) = -y$. Thus, we define \\
\[c_S(||\cdot||) := \sup_{x\dashv_By}\left(\sup_{z\in T_{xy}(S)}||z|| - \inf_{w\in T_{xy}(S)}||w||\right). \] \\
We head now to produce sharp lower and upper bounds for $c_S$, but before this we need a geometric lemma.

\begin{lemma}\label{lemmaball} Let $(V,||\cdot||)$ be a normed plane, and let $x,y \in V$ be such that $x \dashv_B y$. Then $||T_{xy}(z)|| \leq 3$ for any $z \in S$. Moreover, equality is only possible if the plane is rectilinear.
\end{lemma}
\noindent
\textbf{Proof.} For simplicity, assume along the proof that $x,y \in S$, and let $z \in S$. Let $\alpha,\beta \neq 0$ be such that $z = \alpha x + \beta y$ (the other cases are obvious). Since $x \dashv_B y$, we have \\
\[1 = ||\alpha x+ \beta y|| = |\alpha|\left|\left|x+\frac{\beta}{\alpha}y\right|\right| \geq |\alpha|.\]

On the other hand, $1 = ||\alpha x + \beta y|| \geq |\beta| - |\alpha|$. It follows that $|\beta| \leq |\alpha| + 1 \leq 2$. Now, $||T_{xy}(z)|| = ||\alpha x - \beta y|| \leq |\alpha| + |\beta| \leq 3$. \\

Suppose now that there exists a unit vector $z$ such that $||T_{xy}(z)|| = 3$. Writing $z = \alpha x + \beta y$, again we have \\
\[ 2|\beta| \leq ||\alpha x + \beta y|| + ||\alpha x - \beta y|| = 4. \]\\
It follows that $|\beta| = 2$ (the inverse inequality was proved above). Now we have \\
\[ 2|\alpha| = ||\alpha x + \beta y + \alpha x - \beta y|| \geq ||\alpha x - \beta y|| - ||\alpha x + \beta y|| = 2, \]\\
and this yields $|\alpha| = 1$. We may assume that $y + 2x \in S$ (the other cases are completely analogous). Since $x,y \in S$, it follows immediately that the segments $[y(y+2x)]$ and $[(-y)(y+2x)]$ are contained in the unit circle. Therefore $S$ is the parallelogram whose vertices are $\pm y$ and $\pm (y+2x)$.

\begin{flushright} $\square$ \end{flushright}

Notice that we always have $\sup_{z\in S}||T_{xy}(z)|| = \left(\inf_{w \in S}||T_{xy}(w)||\right)^{-1}$. The next corollary follows from this observation. 

\begin{coro}\label{corocs} Let $(V,||\cdot||)$ be a Minkowski plane and fix vectors $x,y \in V$ such that $x \dashv_B y$. Then $\sup_{z\in S}||T_{xy}(z)|| \leq 3$ and $\inf_{w\in S}||T_{xy}(w)|| \geq \frac{1}{3}$. In both cases, equality occurs if and only if the plane is rectilinear.
\end{coro}
\noindent\textbf{Proof.} In view of Lemma 5.1 the proof is straightforward. Notice that a compactness argument shows that the supremum and the infimum are indeed attained for some $z,w \in S$.
\begin{flushright} $\square$ \end{flushright}

\begin{teo}\label{teocs} In any normed plane $(V,||\cdot||)$ we have \\
\[ 0 \leq c_S(||\cdot||) \leq \frac{8}{3}. \] \\
Equality on the left side holds if and only if the norm is derived from an inner product, and equality on the right side holds if and only if the plane is rectilinear.
\end{teo}
\noindent
\textbf{Proof.} The right side follows from Corollary \ref{corocs}. For the left side, if $c_S(||\cdot||) = 0$, then the unit circle is invariant with respect to $T_{xy}$ whenever $x \dashv_B y$. By Theorem \ref{teoroberts} it follows that $x \dashv_R y$ whenever $x\dashv_B y$. This is a characterization of the Euclidean plane (see \cite{alonso}).

\begin{flushright} $\square$ \end{flushright}

We finish by outlining an example where we can calculate the constant $c_S$. The proof is long and very technical, and so we will not present it here.

\begin{figure}
\centering
\includegraphics{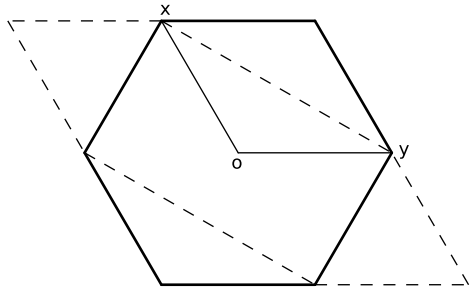}
\caption{$S$ and $T_{xy}(S)$}
\label{fig5e}
\end{figure}

\begin{example}\normalfont If $(V,||\cdot||)$ is a normed plane whose unit circle is an affine regular hexagon, then $c_S(||\cdot||)$ is attained, for example, whenever $y$ is in the direction of one of its sides, $x$ is in the direction of a vertex of this side, and its value equals $\frac{3}{2}$. It is worth mentioning that the constant $c_E$ defined in \cite{bmt} has the same value in regular hexagonal planes. In some sense, this plane is as far from being Euclidean, as Birkhoff and Roberts orthogonal are far from each other. Figure \ref{fig5e} illustrates this ($T_{xy}(S)$ is the dotted polygon).
\end{example}

\bibliography{bibliographyrobertspaper}

\end{document}